\documentclass[reqno,12pt,a4paper]{amsart}

\usepackage{mathrsfs}
\usepackage{amssymb}
\usepackage{amsfonts}
\usepackage{latexsym}
\usepackage{amsthm}

\usepackage{graphicx}
\def\lb{\label}

\newcommand{\er}[1]{\textrm{(\ref{#1})}}

\begin{document}


\renewcommand{\theequation}{\arabic{section}.\arabic{equation}}
\theoremstyle{plain}
\newtheorem{theorem}{\bf Theorem}[section]
\newtheorem{lemma}[theorem]{\bf Lemma}
\newtheorem{corollary}[theorem]{\bf Corollary}
\newtheorem{proposition}[theorem]{\bf Proposition}
\newtheorem{definition}[theorem]{\bf Definition}
\newtheorem{remark}[theorem]{\it Remark}

\def\a{\alpha}  \def\cA{{\mathcal A}}     \def\bA{{\bf A}}  \def\mA{{\mathscr A}}
\def\b{\beta}   \def\cB{{\mathcal B}}     \def\bB{{\bf B}}  \def\mB{{\mathscr B}}
\def\g{\gamma}  \def\cC{{\mathcal C}}     \def\bC{{\bf C}}  \def\mC{{\mathscr C}}
\def\G{\Gamma}  \def\cD{{\mathcal D}}     \def\bD{{\bf D}}  \def\mD{{\mathscr D}}
\def\d{\delta}  \def\cE{{\mathcal E}}     \def\bE{{\bf E}}  \def\mE{{\mathscr E}}
\def\D{\Delta}  \def\cF{{\mathcal F}}     \def\bF{{\bf F}}  \def\mF{{\mathscr F}}
\def\c{\chi}    \def\cG{{\mathcal G}}     \def\bG{{\bf G}}  \def\mG{{\mathscr G}}
\def\z{\zeta}   \def\cH{{\mathcal H}}     \def\bH{{\bf H}}  \def\mH{{\mathscr H}}
\def\e{\eta}    \def\cI{{\mathcal I}}     \def\bI{{\bf I}}  \def\mI{{\mathscr I}}
\def\p{\psi}    \def\cJ{{\mathcal J}}     \def\bJ{{\bf J}}  \def\mJ{{\mathscr J}}
\def\vT{\Theta} \def\cK{{\mathcal K}}     \def\bK{{\bf K}}  \def\mK{{\mathscr K}}
\def\k{\kappa}  \def\cL{{\mathcal L}}     \def\bL{{\bf L}}  \def\mL{{\mathscr L}}
\def\l{\lambda} \def\cM{{\mathcal M}}     \def\bM{{\bf M}}  \def\mM{{\mathscr M}}
\def\L{\Lambda} \def\cN{{\mathcal N}}     \def\bN{{\bf N}}  \def\mN{{\mathscr N}}
\def\m{\mu}     \def\cO{{\mathcal O}}     \def\bO{{\bf O}}  \def\mO{{\mathscr O}}
\def\n{\nu}     \def\cP{{\mathcal P}}     \def\bP{{\bf P}}  \def\mP{{\mathscr P}}
\def\r{\rho}    \def\cQ{{\mathcal Q}}     \def\bQ{{\bf Q}}  \def\mQ{{\mathscr Q}}
\def\s{\sigma}  \def\cR{{\mathcal R}}     \def\bR{{\bf R}}  \def\mR{{\mathscr R}}
\def\S{\Sigma}  \def\cS{{\mathcal S}}     \def\bS{{\bf S}}  \def\mS{{\mathscr S}}
\def\t{\tau}    \def\cT{{\mathcal T}}     \def\bT{{\bf T}}  \def\mT{{\mathscr T}}
\def\f{\phi}    \def\cU{{\mathcal U}}     \def\bU{{\bf U}}  \def\mU{{\mathscr U}}
\def\F{\Phi}    \def\cV{{\mathcal V}}     \def\bV{{\bf V}}  \def\mV{{\mathscr V}}
\def\P{\Psi}    \def\cW{{\mathcal W}}     \def\bW{{\bf W}}  \def\mW{{\mathscr W}}
\def\o{\omega}  \def\cX{{\mathcal X}}     \def\bX{{\bf X}}  \def\mX{{\mathscr X}}
\def\x{\xi}     \def\cY{{\mathcal Y}}     \def\bY{{\bf Y}}  \def\mY{{\mathscr Y}}
\def\X{\Xi}     \def\cZ{{\mathcal Z}}     \def\bZ{{\bf Z}}  \def\mZ{{\mathscr Z}}
\def\O{\Omega}

\newcommand{\mc}{\mathscr {c}}

\newcommand{\gA}{\mathfrak{A}}          \newcommand{\ga}{\mathfrak{a}}
\newcommand{\gB}{\mathfrak{B}}          \newcommand{\gb}{\mathfrak{b}}
\newcommand{\gC}{\mathfrak{C}}          \newcommand{\gc}{\mathfrak{c}}
\newcommand{\gD}{\mathfrak{D}}          \newcommand{\gd}{\mathfrak{d}}
\newcommand{\gE}{\mathfrak{E}}
\newcommand{\gF}{\mathfrak{F}}           \newcommand{\gf}{\mathfrak{f}}
\newcommand{\gG}{\mathfrak{G}}           
\newcommand{\gH}{\mathfrak{H}}           \newcommand{\gh}{\mathfrak{h}}
\newcommand{\gI}{\mathfrak{I}}           \newcommand{\gi}{\mathfrak{i}}
\newcommand{\gJ}{\mathfrak{J}}           \newcommand{\gj}{\mathfrak{j}}
\newcommand{\gK}{\mathfrak{K}}            \newcommand{\gk}{\mathfrak{k}}
\newcommand{\gL}{\mathfrak{L}}            \newcommand{\gl}{\mathfrak{l}}
\newcommand{\gM}{\mathfrak{M}}            \newcommand{\gm}{\mathfrak{m}}
\newcommand{\gN}{\mathfrak{N}}            \newcommand{\gn}{\mathfrak{n}}
\newcommand{\gO}{\mathfrak{O}}
\newcommand{\gP}{\mathfrak{P}}             \newcommand{\gp}{\mathfrak{p}}
\newcommand{\gQ}{\mathfrak{Q}}             \newcommand{\gq}{\mathfrak{q}}
\newcommand{\gR}{\mathfrak{R}}             \newcommand{\gr}{\mathfrak{r}}
\newcommand{\gS}{\mathfrak{S}}              \newcommand{\gs}{\mathfrak{s}}
\newcommand{\gT}{\mathfrak{T}}             \newcommand{\gt}{\mathfrak{t}}
\newcommand{\gU}{\mathfrak{U}}             \newcommand{\gu}{\mathfrak{u}}
\newcommand{\gV}{\mathfrak{V}}             \newcommand{\gv}{\mathfrak{v}}
\newcommand{\gW}{\mathfrak{W}}             \newcommand{\gw}{\mathfrak{w}}
\newcommand{\gX}{\mathfrak{X}}               \newcommand{\gx}{\mathfrak{x}}
\newcommand{\gY}{\mathfrak{Y}}              \newcommand{\gy}{\mathfrak{y}}
\newcommand{\gZ}{\mathfrak{Z}}             \newcommand{\gz}{\mathfrak{z}}

\def\ve{\varepsilon}   \def\vt{\vartheta}    \def\vp{\varphi}
\def\vk{\varkappa}

\def\A{{\mathbb A}} \def\B{{\mathbb B}} \def\C{{\mathbb C}}
\def\dD{{\mathbb D}} \def\E{{\mathbb E}} \def\dF{{\mathbb F}}
\def\dG{{\mathbb G}} \def\H{{\mathbb H}}\def\I{{\mathbb I}}
\def\J{{\mathbb J}} \def\K{{\mathbb K}} \def\dL{{\mathbb L}}
\def\M{{\mathbb M}} \def\N{{\mathbb N}} \def\dO{{\mathbb O}}
\def\dP{{\mathbb P}} \def\R{{\mathbb R}}\def\S{{\mathbb S}}
\def\T{{\mathbb T}} \def\U{{\mathbb U}} \def\V{{\mathbb V}}
\def\W{{\mathbb W}} \def\X{{\mathbb X}} \def\Y{{\mathbb Y}}
\def\Z{{\mathbb Z}}


\def\la{\leftarrow}              \def\ra{\rightarrow}  \def\Ra{\Rightarrow}
\def\ua{\uparrow}                \def\da{\downarrow}
\def\lra{\leftrightarrow}        \def\Lra{\Leftrightarrow}


\def\lt{\biggl}                  \def\rt{\biggr}
\def\ol{\overline}               \def\wt{\widetilde}
\def\no{\noindent}


\let\ge\geqslant                 \let\le\leqslant
\def\lan{\langle}                \def\ran{\rangle}
\def\/{\over}                    \def\iy{\infty}
\def\sm{\setminus}               \def\es{\emptyset}
\def\ss{\subset}                 \def\ts{\times}
\def\pa{\partial}                \def\os{\oplus}
\def\om{\ominus}                 \def\ev{\equiv}
\def\iint{\int\!\!\!\int}        \def\iintt{\mathop{\int\!\!\int\!\!\dots\!\!\int}\limits}
\def\el2{\ell^{\,2}}             \def\1{1\!\!1}
\def\sh{\sharp}
\def\wh{\widehat}
\def\bs{\backslash}
\def\intl{\int\limits}

\def\na{\mathop{\mathrm{\nabla}}\nolimits}
\def\sh{\mathop{\mathrm{sh}}\nolimits}
\def\ch{\mathop{\mathrm{ch}}\nolimits}
\def\where{\mathop{\mathrm{where}}\nolimits}
\def\all{\mathop{\mathrm{all}}\nolimits}
\def\as{\mathop{\mathrm{as}}\nolimits}
\def\Area{\mathop{\mathrm{Area}}\nolimits}
\def\arg{\mathop{\mathrm{arg}}\nolimits}
\def\const{\mathop{\mathrm{const}}\nolimits}
\def\det{\mathop{\mathrm{det}}\nolimits}
\def\diag{\mathop{\mathrm{diag}}\nolimits}
\def\diam{\mathop{\mathrm{diam}}\nolimits}
\def\dim{\mathop{\mathrm{dim}}\nolimits}
\def\dist{\mathop{\mathrm{dist}}\nolimits}
\def\Im{\mathop{\mathrm{Im}}\nolimits}
\def\Iso{\mathop{\mathrm{Iso}}\nolimits}
\def\Ker{\mathop{\mathrm{Ker}}\nolimits}
\def\Lip{\mathop{\mathrm{Lip}}\nolimits}
\def\rank{\mathop{\mathrm{rank}}\limits}
\def\Ran{\mathop{\mathrm{Ran}}\nolimits}
\def\Re{\mathop{\mathrm{Re}}\nolimits}
\def\Res{\mathop{\mathrm{Res}}\nolimits}
\def\res{\mathop{\mathrm{res}}\limits}
\def\sign{\mathop{\mathrm{sign}}\nolimits}
\def\span{\mathop{\mathrm{span}}\nolimits}
\def\supp{\mathop{\mathrm{supp}}\nolimits}
\def\Tr{\mathop{\mathrm{Tr}}\nolimits}
\def\BBox{\hspace{1mm}\vrule height6pt width5.5pt depth0pt \hspace{6pt}}


\newcommand\nh[2]{\widehat{#1}\vphantom{#1}^{(#2)}}
\def\dia{\diamond}

\def\Oplus{\bigoplus\nolimits}



\def\qqq{\qquad}
\def\qq{\quad}
\let\ge\geqslant
\let\le\leqslant
\let\geq\geqslant
\let\leq\leqslant
\newcommand{\ca}{\begin{cases}}
\newcommand{\ac}{\end{cases}}
\newcommand{\ma}{\begin{pmatrix}}
\newcommand{\am}{\end{pmatrix}}
\renewcommand{\[}{\begin{equation}}
\renewcommand{\]}{\end{equation}}
\def\eq{\begin{equation}}
\def\qe{\end{equation}}
\def\[{\begin{equation}}
\def\bu{\bullet}

\title[{Estimates of 1D resonances in terms of potentials }]
{Estimates of 1D resonances in terms of potentials}

\date{\today}

\author[Evgeny Korotyaev]{Evgeny Korotyaev}
\address{Mathematical Physics Department, Faculty of Physics, Ulianovskaya 2,
St. Petersburg State University, St. Petersburg, 198904, Russia,
 \ korotyaev@gmail.com
}


\begin{abstract}
\no We discuss resonances for Schr\"odinger operators with compactly supported
potentials on the line and the half-line. We estimate the  sum of the negative
 power of
all resonances and eigenvalues in terms of the norm of the potential
and the diameter of its support. The proof is based on harmonic
analysis and Carleson measures arguments.
\end{abstract}

\maketitle

{\it Dedicated to Lennart Carleson, on the occasion of his 85th birthday}

\section {Introduction and main results}
\setcounter{equation}{0}

In this paper we will present a global estimate of resonances in terms of the
potential for Schr\"odinger operators   $H=H_0+q$, where $H_0$ is one of the
following:
$$
\begin{aligned}
&Case \ 1: \qqq -{d^2/dx^2}\qq {\rm in}\qq L^2(\R ).\\
&Case \ 2: \qqq -{d^2/dx^2}\qq {\rm in}\qq L^2(\R_+),\qq
{\rm with} \ f(0)=0 \ {\rm boundary\ conditions}.\\
&Case \ 3: \qqq -{d^2/dx^2}\qq {\rm in}\qq L^2(\R_+),\qq
{\rm with}\ \ f'(0)=0 \ {\rm boundary\ conditions}.
\end{aligned}
$$
We assume that $q$ is real, integrable and has a compact support. It
is well known that the spectrum of $H$ consists of an absolutely
continuous part $[0,\iy)$ and a finite number of simple negative
eigenvalues $E_1<\dots <E_m<0$, see well-known papers \cite{DT79},
\cite{Me85} and the book \cite{Ma86} about inverse scattering. The
Schr\"odinger equation
\[
\lb{e1a}
  -f''+q(x)f=k^2f,\ \ \  \ \ \ k\in \C \sm \{0\},
\]
has unique solutions $\p_{\pm}(x,k)$  such that $\p_+(x,k )=e^{ixk }$  for large
positive $x$ and $\p_-(x,k)=e^{-ik x}$ for large negative $x$.
Outside the support of $q$ any solutions of \er{e1a} have to be combinations
of $e^{\pm ik x}$. The functions $\p_\pm(x,\cdot), \p_\pm'(x,\cdot)$ for all $x\in \R$
are entire. We define the Wronskian $w$ for Case 1 by
\[
\begin{aligned}
\lb{H1}
w(k )=\{\p_-(\cdot,k), \p_+(\cdot,k )\},\qqq
\end{aligned}
\]
where  $\{f, g\}=fg'-f'g$. In  Case 2 the Jost function is defined
as $\p_+(0,\cdot)$  and  in Case 3 the Jost function is defined as
$\p_+'(0,\cdot)$. Let $F$ be one of the functions $w, \p_+(0,k )$ or
$\p_+'(0,k )$. Recall that the function $F$ is entire. One  has
exactly $m$ simple zeros $k_{1}=i|E_1|^{1/2},\dots
k_{m}=i|E_{m}|^{1/2}$ in the upper half-plane $\C_+$ and for $q\ne
0$ an infinite number of zeros $(k_{n})_{m+1}^{\iy }$ in the lower
half-plane  $\ol\C_-$ labeled  by
$$
0\leq |k_{m+1}|\leq |k_{m+2}|\leq \dots \qqq where \qqq 0<|k_{m+2}|,
$$
see \cite{F97},  \cite{S00}, \cite{Z87} and \cite{K05}.
Here it is possible that $k_{m+1}=0$ for some potential, but $0<|k_{m+2}|$
for any potential. By definition, a zero $k_n\in \ol\C_-$ of $F$ is called a resonance
of $H$.  The multiplicity of the resonance is the multiplicity
of the corresponding zero of $F$.
Of course, the energies are given by $k^2$, but since $k$ is the natural parameter,
we will abuse the terminology.

 There are only few estimates of resonances.
 We denote the number of zeros of function $f$
having modulus $\le r$ by $\cN (r,f)$, each zero being counted
according to its multiplicity. Firstly, Zworski \cite{Z87}
determines the asymptotics of the counting function for resonances :
\[
\lb{zi} \cN (r,w)={2 r\/\pi}(\g+o(1))\qqq \as \qq r\to \iy,
\]
where  $[0,\g]$ is the convex hull of the support of $q$.

Secondly, let $\mN_t$ denote the number of
resonances and eigenvalues in the half-plane $\{\Im k>t\}, t<0$.
Then there is a constant $\cC_q$  (see Theorem 3.11 in \cite{H99})
such that the following estimate holds true:
\[
\lb{Hi}
\mN_t\le
\cC_q\rt(1+\int\int_{\R^2}e^{4|t||x-y|}|q(x)||q(x)|dxdy\rt),
\]
where $\cC_q$ is some constant depending on $\|q\|$, but not on $t$.
Unfortunately, \er{Hi} is not sharp, since  $\cC_q$ is unknown.

Thirdly, for the Case 2 (at $\g=1$) there are simple estimates from \cite{K04}
\[
\lb{kn04} |k_n|e^{-2|\Im k_n|}\le \|q\|e^{\|q\|},\qqq \qqq  {\rm
where} \ \ \|q\|=\int_\R |q(x)|dx,
\]
 for  any $k_n\in\C_-$. Note that this estimate yields the
will-known logarithmic curve for forbidden domain. If, in addition,
$q'$ is integrable, then \er{kn04} will be sharper, and so on.

Define the constant $Q$ by
\[
Q=\max\{\|q\|,\|q\|_1\},\qqq \|q\|=\int_\R |q(t)|dt,\qqq
\|q\|_1=\int_\R |tq(t)|dt.
\]

We present theorem about new estimates of counting functions.

\begin{theorem}
\lb{TrJS} Let $H=H_0+q$, where $q$ is integrable and has a compact
support.

In Case 1, let $\supp q\ss [0, \g]$ but in no smaller interval.

In Cases 2 and 3, let $\g=\sup (\supp(q))$.

 Let $r>0$ and $r_1=r+{1\/2}$. Then the following estimates
 hold true:
 \[
\lb{q2} \cN (r,w)\le {1\/\log 2}\rt({4r_1\g\/ \pi }+\log (1+{4r_1})+
{9Q\/1+4r_1}\rt),
\]
\[
\lb{q3}
{\cN (r,\p_+(0,\cdot))+\cN (r,\p_+'(0,\cdot))\/2}\le
{1\/\log 2}\rt({4r_1\g\/ \pi }+\log (1+{4r_1})+{9\|q\|\max\{1,\g\}\/1+4r_1}\rt).
\]
\end{theorem}

 The proof is based on the Jensen formula and
standard estimates of the fundamental solutions. The RHS in \er{q2}
has asymptotics ${2\/\log 2}{2 r\/\pi}(\g+o(1))$ as $r\to \iy$. If
we compare this asymptotics  with \er{zi}, then we obtain the
coefficient ${2\/\log 2}$. It means that the estimate \er{q2} is
sufficiently sharp.

We present our main result.

\begin{theorem}
\lb{T1}
Let $H=H_0+q$, where $q$ is integrable and has a compact support.

In Case 1, let $\supp q\ss [0, \g]$ but in no smaller interval.

In Cases 2 and 3, let $\g=\sup (\supp(q))$.

Then for any $p>1$ the following estimates hold true:
\[
\begin{aligned}
\lb{3} \sum_{\pm \Im k_n\le 0} {1\/| k_n-2i|^{p}}\le
CY_p\rt(1+{\g\/ \pi }+\cQ\rt),\\
\end{aligned}
\]
where $C\le 2^5$ is an absolute  constant,
$Y_p=\sqrt \pi{\G({p-1\/2})/\G({p\/2})}$ and
\[
\cQ=\ca  \max\{\|q\|,\|q\|_1\} & Case \ 1\\
           2\|q\|\max \{1, \g\},                    & Case\  2,3
\ac.
\]
\end{theorem}

{\bf Remark.} 1) The function $Y_p, p>1$  is strongly monotonic and
convex on $(1,\iy)$, since $Y_p'<0, Y_p''>0$, and satisfies (see
more in Lemma \ref{Tp})
\[
\lb{Yp} Y_2=\pi ,\qqq
Y_p=\ca {1\/\sqrt p}(\sqrt {2\pi}+O(1/p))  & \as \qq p\to \iy\\
             {1\/p-1}(2+o(1)) & \as \qq p\to  1\ac.
\]
Thus we can control the RHS of \er{3} at $p\to 1$ and for large
$p\to \iy$. Note we take $p>1$, since  the asymptotics \er{zi}
implies the simple fact $\sum_{k_n\ne 0} {1\/|k_n|}=\iy$, see p. 17
in \cite{L93}.

2) The RHS of \er{3} depends on 3 crucial parameters: $p>1$, the
diameter of the support of the potential and the magnitude $\|q\|$
of the potential. In Cases 1 and 3 we can not remove 1 in the RHS of
\er{3}. In Case 2 probably the number 1 should be  absent in the RHS
of \er{3}. In order to explain this we need to add that at $q=0$
there is a resonance in the Cases 1 and 3, but there is no a
resonance in the Case 2.

3) In Case 1 the proof of \er{3} is based on analysis of the
function $w$. We use harmonic  analysis  and the Carleson Theorem (
Theorems 1.56 and 2.3.9 in \cite{G81}) about Carleson measure. The
proof for Cases 2 and 3 is a simple corollary of Case 1.

4)   $C$ is the constant from Carleson's Theorem ( \cite{C58},
\cite{C62} and  see Theorems 1.56 and   2.3.9, \cite{G81}), see also
\er{Ce}.

5) In fact, the estimates \er{q2},\er{3} give a new global property
of resonance stability.

6) It is well-known that the Jost function is a Fredholm
determinant, see \cite{JP51}. In Case 1 we consider the Fredholm
determinant $D(k)=\det (I+q(H_0-k^2)^{-1}, k\in \C_+$. The function
$D$ is analytic in the upper half-plane (see e.g. \cite{F97}) and
satisfies the well-known identity $D(k)={w(k)\/2ik}$ for all $k$.
Thus Theorem \ref{T1} describes the zeros of  the determinants also.

\medskip

Resonances for the multidimensional case were studied by Melrose,
Sj\"ostrand, and Zworski and other, see \cite{M83}, \cite{Z89},
\cite{SZ91}) and references therein. We discuss the one dimensional
case.  A lot of papers is devoted to resonances for the 1D
Schr\"odinger operator, see Froese \cite{F97}, Simon \cite{S00},
Zworski \cite{Z87}, Korotyaev \cite{K11} and references therein.
Different properties of resonances were determined in \cite{H99},
\cite{S00}, \cite{Z87}, and  \cite{K04}, \cite{K05}, \cite{K11}.
Korotyaev solved the inverse problem for resonances for the
Schr\"odinger operator  with a compactly  supported potential on the
real line \cite{K05} and the half-line \cite{K04}: (i) the
characterization  of $S$-matrix in terms of resonances, (ii) a
recovering of the potential from the resonances, (iii)
 the potential is uniquely determined by the resonances (about uniqueness
  see also \cite{Z02}, \cite{BKW03}).

The "local resonance" stability problem was considered in
\cite{K04s} for  Case 2. Roughly speaking, if $(k_n)_1^\iy$ is a
sequence of eigenvalues and resonances of the Schr\"odinger operator
with some compactly supported potential $q$ and $\sum _{n\ge
1}n^{2\ve}|k_n-k_{n}^\bullet|^2<\iy$ for some  sequence
$(k_n^\bullet)_1^\iy$ and $\ve>1$, then $(k_n^\bullet)_1^\iy$ is a
sequence of eigenvalues and resonances of a Schr\"odinger operator
for some unique real compactly supported potential $ q_\bullet$.

Finally we note that some  stability results (quite different type,
so-called stability for finite data) were obtained by
Marlettta-Shterenberg-Weikard \cite{MSW10}.

A lot of papers are devoted to estimates of eigenvalues in terms of
the  potential. We recall Lieb-Thirring type inequalities (\cite{LT76}) from \cite{DFLP06} given by:

Let $V$ be a non-negative, unbounded potential, such that
the Schr\"odinger operator $-\D+V$ in $L^2(\R^d)$
has an unbounded sequence of discrete  eigenvalues $\l_1\le \l_2\le...$. Then
\[
\lb{LT}
 \sum_{n\ge 1} |\l_n|^{-{p\/2}}\le   {C_p(d)\/(4p)^{d\/2}} \int_{\R^d} V(x)^{{d-p\/2}}dx,\qqq
 where \ C_p(d)={\G({p-d\/2})\/\G({p\/2})},
\]
 for all $ p>d$. The inequality is derived from inequalities by Golden
 \cite{G65} and Thompson \cite{T65}. It is interesting that \er{3}
 and \er{LT} (at $d=1$) have the same constant $C_p(1)=Y_p$.

In our case \er{3}  the number of resonances is infinite and we
consider also the  power $p>1$. Furthermore, some resonances and
eigenvalues can be close to zero even for small potentials. For this
reason,   we sum $|k_n-2i|^{-{p}}$ and $|k_n+2i|^{-{p}}$.  Thus,
roughly speaking,  \er{3} is a Lieb-Thirring type inequality for the
resonances \cite{L13}.




\section {Proof}
\setcounter{equation}{0}

\subsection{Estimates for entire functions.} An entire
function $f(z)$ is said to be of exponential type if there is a
constant $\b$ such that $|f(z)|\leq\const e^{\b |z|}$ everywhere.
The infimum of the set of $\b$ for which such inequality holds is
called the type of $f$.

\no {\bf Definition.}
{\it Let $\cE_\g, \g >0$ denote the set
 of exponential type  functions $f$, which satisfy
\[
\lb{D1} |f(k)|\ge 2|k|\qqq   \qqq \forall  \ k\in \R,
\]
\[
\lb{D2}
\begin{aligned}
&|f(k)-2ik+f_0|\le {Q^2\/|k|_1}e^{ \g (|\Im k|-\Im
k)+{Q\/|k|_1}},\qqq
\qq  \forall \ k\in \C,\\
& Q> 0,\qq  f_0\in \R,\qqq \qqq |f_0|\le Q,
\end{aligned}
\]
where the constants $Q=Q(f), f_0$ depend on $f$ and
$|k|_1=\max\{1,|k|\}$.}

In the proof of Theorem \ref{T1} we will need some properties of zeros of
$f\in \cE_\g $ in terms of the Carleson  measure. Recall that {\it a
positive Borel measure $M$ defined in $\C_-$ is called a
Carleson measure if there is
 a constant $C_M$ such that for all $(r,t)\in \R_+\ts\R$
\[
\lb{1.31}
   M(D_-(t,r))\leq C_Mr,\ \ \ {\rm where}\ \ \
   \ D_-(t,r)\ev\{z\in \C_-: |z-t|<r \},
\]
here $C_M$ is the Carleson  constant independent of $(t,r)$.}

For an entire function $f$ with zeroes $k_n, n\ge 1$
  we define an associated measure by
\[
\lb{MO}
 d\O (k,f)=\sum_{\Im k_n\le 0} \d (k-k_n+i)du dv,\qqq k=u+iv\in \C_-.
\]
We denote the number of zeros of function $f$ having modulus  $<r$ by $\cN (r,f)$,
each zero being counted according to its multiplicity.
In order to prove Theorem \ref{T1} we need

  \begin{theorem}
\lb{TEf}
 Let $f\in \cE_\g $ for some $\g >0$. Then

 i) The number of zeros of $f$ in the disk $\{|k|<r\},r>0, r_1=r+{1\/2}$
satisfies
\[
\lb{L1} \cN (r,f)\le {1\/\log 2}\rt({4r_1\g\/ \pi }+\log (1+{4r_1})+
{|f_0|+8Q\/1+4r_1}\rt).
\]

ii) $\O (\cdot,f)$ is a Carleson measure and for all
$(r,t)\in \R_+\ts \R $ satisfies
\[
\begin{aligned}
\lb{L2}
\O (D_-(t,r),f)\le\cN (r,f(t+\cdot))\le C(f)r ,\\
C(f)={12\/\log 4}\rt({\g \/\pi }+1+{|f_0| +3Q\/4}\rt).
\end{aligned}
\]
iii) Let $k_n, n\ge 1$ be all zeros of $f$. Then for each $p>1$ the following
estimate holds true:
\[
\lb{L3}
\sum_{\pm \Im k_n\le 0} {1\/|k_n-2i|^{p}}
\le   C Y_pC(f),
\]
where $C\le 2^5$ is the absolute  constant and  $Y_p$ is given by
\[
\lb{p1}
  Y_p=\int_\R
(1+x^2)^{-{p\/2}}dx=\sqrt \pi{\G({p-1\/2})\/\G({p\/2})} \qqq \forall \qq p>1.
\]

\end{theorem}
{\bf Proof.} i) Recall the Jensen formula (see p. 2 in \cite{Koo88})
for an entire function $F$:
\[
\lb{PJ}
\log |F(0)|+\int _0^r{\cN (t,F)\/ t}dt= {1\/ 2\pi }\int
_0^{2\pi}\log |F(re^{i\f})|d\f,
\]
 Rewrite $f$ in the form
 \[
 \lb{ef1}
 f=2ik-f_0+f_*, \qqq where \qq |f_*(k)|\le
{Q^2\/|k|_1}e^{{Q\/|k|_1}}e^{2\g v_-},\qqq v_-(k)=
{|\Im k|-\Im k\/2},
\]
for all $k\in \C$.
Consider the case $t\in \R\sm [-{1\/2}, {1\/2}]$. Let $z=t+k, k=re^{i\f}$.
Then \er{D2} and \er{ef1} imply
\[
\begin{aligned}
|f(t+k)|=|2it+2ik-f_0+f_*(z)|\le 2|t|\rt(1+{r\/|t|}\rt)\rt(1+{|f_0|+|f_*(z)|\/2(|t|+r)}\rt)\\
\le 2|t|\rt(1+{r\/|t|}\rt) e^{2\g v_-+{Q\/|z|_1}}\rt(1+{|f_0|+Q^2|z|_1^{-1}\/2(|t|+r)}\rt)\\
\le 2|t|\rt(1+{r\/|t|}\rt) e^{2\g v_-+{Q\/|z|_1}}
\rt(1+{|f_0|\/2(|t|+r)}\rt)  \rt(1+{Q^2\/2(|t|+r)|z|_1}\rt).
\end{aligned}
\]
Let $\cN (s)=\cN (s,f(t+\cdot))$.
Substituting this into \er{PJ} we obtain for the function $F(k)=f(t+k),
k=re^{i\f}, z=t+re^{i\f}$:
\[
\begin{aligned}
\lb{F1}
\log |f(t)|+\int_0^r{\cN (s)\/s}ds={1\/2\pi}\int_0^{2\pi}\log |f(z)|d\f\\
\le \log 2|t|+\log \rt(1+{r\/|t|}\rt)+\log \rt(1+{|f_0|\/2(|t|+r)}\rt)+{1\/2\pi}\int_0^{2\pi}2\g v_-d\f+X(t,r),\\
where \qq X(t,r)={1\/2\pi}\int_0^{2\pi}\rt[{Q\/|z|_1}+\log  \rt(1+{Q^2\/2(|t|+r)|z|_1}\rt) \rt]d\f,
\end{aligned}
\]
where $2v_-=r(|\sin \f| -\sin\f)$ and the simple integration yields
\[
\lb{F12}
\begin{aligned}
{1\/2\pi}\int_0^{2\pi}2\g v_-d\f=
{r\g\/2\pi }\int_0^{2\pi}\rt(|\sin \f|-\sin \f\rt)d\f={2r\g\/\pi}.
\end{aligned}
\]
Substituting the estimate ${|f(t)|\/2|t|}\ge 1$  into \er{F1}
together  with the simple one
$$
\int _0^r{\cN (s)ds\/s}\geq \cN ({r/2})\int _{r\/2}^r{ds\/s}=
\cN ({r/2})\log 2,
$$
we obtain
\[
\begin{aligned}
\lb{F3}
\cN ({r/2})\log 2\le {2r\g\/ \pi }+
\log \rt(1+{r\/|t|}\rt)+{|f_0|\/2(|t|+r)}+X(t,r).
\end{aligned}
\]

We show \er{L1}. Let $t={1\/2}$ and $r_1=r+{1\/2}$ for any $r>0$.
Then
$$
\{|k|<r\}\ss \{|k-{1\/2}|<r_1\},\qqq |z|_1=|(1/2)+k|_1\ge {1\/2}(1+|k|-(1/2))=
{1\/4}(1+2|k|)
$$
and \er{F3} give
\[
\lb{esr}
\begin{aligned}
\cN (r,f)\le \cN (r_1,f({1\/2}+\cdot))\le
{4r_1\g\/ \pi }+\log (1+4r_1)+{|f_0|\/(1+4r_1)}+X(t,r),\\
where \qq X(t,r)\le {1\/2\pi}\int_0^{2\pi}\rt[{4Q\/(1+4r_1)}+\log  \rt(1+{4Q^2\/(1+4r_1)^2}\rt) \rt]d\f\le {8Q\/(1+4r_1)},
\end{aligned}
\]
which yields  \er{L1}.

ii) 
Let $r\le 1, t\in \R$. Then by the construction of $\O (\cdot,f)$, we obtain
$\O (D_-(t,r),f) =0$.

  In order to show \er{L2} for thew case  $r> 1, t\in \R$ we need to consider two cases:

Firstly, let  $r>1, t\in \R\sm [-{1\/2}, {1\/2}]$. Then due to \er{F3} and \er{F1},
the measure $\O (\cdot,f)$ satisfies
\[
\begin{aligned}
\lb{F4}
& \O (D_-(t,r),f)\le\cN (r,f(t+\cdot))\le
{1\/\log 2}\rt(   {4r\g\/ \pi }+4r+|f_0| +Q+\log (1+Q)^2\rt)\le C_1r,\\
& \qqq \qqq  C_1={4\/\log 2}\rt({\g\/\pi }+1+{|f_0| +3Q\/4}\rt).
\end{aligned}
\]

Secondly, let  $r>1, t\in [0, {1\/2}]$. The proof for the case $t\in [-{1\/2},0]$
is similar. For two disks $D(t,r)$ and $D({1\/2},r_1), r_1=r+{1\/2}$ we have
\[
\begin{aligned}
\lb{F5}
D(t,r)\cap \{\Im k\le -1\}\ss D\big({1\/2},r_1\big)\cap \{\Im k\le -1\},
\\
r\le r_1=r+{1\/2}\le {3r\/2}.
\end{aligned}
\]
Then due to \er{L1}, the measure $\O (\cdot,f)$ satisfies
\[
\begin{aligned}
\lb{F4a}
\O (D_-(t,r),f)\le \O (D_-(1/2,r_1),f)\le \cN (r_1,f({\small{1\/2}}+\cdot))\le
C_1r_1\le {3\/2}C_1r.
\end{aligned}
\]
Thus $\O (\cdot,f)$ is a Carleson measure with the Carleson constant $C(f)={3\/2}C_1$, where $C_1$ is given by \er{F4}.

iii)  Consider the case $\Im k_n\le 0$. The proof for the case $\Im k_n> 0$
is similar, even simpler.
 In order to show \er{L3} we recall the Carleson result (see p.
63, Theorem 3.9, \cite{G81}):

Let $F$ be analytic on $\C_-$. For $0<p<\iy$  we say $F\in
\mH_p=\mH_p(\C_-)$ if
$$
\sup_{y<0}\int_\R|F(x+iy)|^pdx=\|F\|_{\mH_p}^p<\iy
$$
Note that the definition of the Hardy space $\mH_p$ involve all
$y<0$, instead of small only value of $y$, like say, $y\in (-1,0)$.
We define the Hardy space $\mH_p$ for the case $\C_-$, since below
we work with functions on $\C_-$.

{\it If $M$ is a Carleson measure and satisfies \er{1.31}, then the
following estimate holds:
\[
\lb{Ce}
  \int_{\C_-} |F|^pdM\leq C C_M \|F\|_{\mH_p}^p\qqq
  \forall \qq F\in \mH_p, \ p\in (0,\iy ),
\]
where  $C_M$ is the so-called Carleson constant from \er{1.31} and  $C\le 2^5$ is an
absolute constant.}

 In order to prove \er{L3} we  take the functions $F(k)={1\/k-i}$.
Then estimates \er{Ce}, \er{L2}  yield
\[
\lb{3.19} \sum_{\Im k_n\le 0} {1\/|k_n-2i|^p}=\int_{\C_-}
|F(\l)|^pd\O \leq C C(f) \|F\|_{\mH_p}^p, \qq p\in (1,\iy),
\]
$C$ is an absolute constant and $C(f)$ is defined in \er{L2}. Here
we have the simple identity
\[
\lb{3.19a} \|F\|_{\mH_p}^p=\int_\R {dt\/|t-i|^p}=\int_\R
{dt\/(t^2+1)^{p\/2}}=Y_p.
\]
The function $Y$ is studied in Lemma \ref{Tp}. Combine \er{3.19} and
\er{3.19a} we obtain \er{L3}.
\hfill  \BBox

\

\subsection{Estimates of resonances.} We consider Case 1 and without loss of
generality let $\supp q\ss [0,\g]$. The Schr\"odinger equation
\[
\lb{e1}
  -f''+q(x)f=k^2f,\ \ \  \ \ \ k\in \C \sm \{0\},
\]
has unique solutions $\p_{\pm}(x,k)$  such that $\p_+(x,k )=e^{ixk
},\ \ x\ge \g$ and $\p_-(x,k)=e^{-ik x},\ \ x\leq 0$. Outside of the
support of $q$ any solutions of \er{e1} have to be linear
combinations of $e^{\pm ik x}$. Functions $\p_\pm(x,\cdot),
\p_\pm'(x,\cdot),  x\in \R$ are entire. We define the functions $a,
w, s$ by
\[
\begin{aligned}
\lb{H1a}
w(k)=2ika(k)=\{\p_-(\cdot,k), \p_+(\cdot,k )\},\qqq
  \ s(k)=\{\p_+(\cdot,k ),\p_-(\cdot,-k )\},
\end{aligned}
\]
and $\{f, g\}=fg'-f'g$ denotes the Wronskian.  The functions $w(k),
s(k)$ are entire and the following asymptotic
estimates hold true:
\[
\lb{a1}
 a(k )=1+O(k ^{-1})\qqq {\rm as} \ \ |k |\to \iy ,\ \  k\in\ol\C_+.
\]
The scattering matrix for the operators $H, H_0=-{d^2\/dx^2}$ is given by
\[
S(k)=\ma a(k )^{-1}& r_-(k )\\ r_+(k )&a(k )^{-1}\am ,
\ \ \ \  r_{\pm}={s(\mp k)\/ w(k )}, \ \ \ \  k\in\R,
\lb{1.7}
\]
where $1/a$ is the transmission coefficient and $r_{\pm}$ are the reflection
coefficients. The matrix $S(k), k\in \R$ is unitary, which yields
\[
\lb{ws}
|w(k)|^2=4k^2+|s(k)|^2,\ \  \forall \ k\in \R.
\]
The solution $ \p_+$ of \er{e1}  satisfies the following equation
\[
\lb{3.1}
 \p_+(x,k)=e^{ixk}-\int _x^\g{\sin [k(x-t)]\/k}q(t)
\p_+(t,k)dt \qqq \forall\  (x,k)\in [0, \g]\ts \C.
\]
It is well known that equation \er{3.1} has a unique solution.
 Due to \er{3.1} the function $y(x,k)=e^{-ikx}\p_+(x,k)$ satisfies
the integral equation
\[
\lb{y1}
y(x,k)=1+\int _x^\g G(t-x,k)q(t)y(t,k)dt \qqq  G(t,k)={\sin kt\/k}e^{ikt},
\]
for all $(x,k)\in [0,\g]\ts \C$. We have the standard iterations
\[
\lb{y2}
y(x,k)=1+\sum_{n\ge 1}y_n(x,k),\qqq y_n(x,k)=
\int _x^\g G(t-x,k)q(t)y_{n-1}(t,k)dt,
\qq y_0=1.
\]
We need some properties of the functions introduced above.

\begin{lemma}
\label{Tw}
 Let $q\in L^1(\R)$ and  $\supp q\ss [0,\g]$. Then the functions
 $w, s$ and $\p _{+}(x,\cdot),\p _{+}'(x,\cdot)$, $x\in\R$ are entire and we have
\[
\lb{yn}
|y_n(x,k)|\le {h^n\/n!}e^{2(\g-x)v_-}, \qqq v_-=\ca 0,&  k\in \ol\C_+\\
|\Im k|,  &  k\in \C_-\ac,
\]
for any $n\ge 1, (x,k)\in [0,\g]\ts \C$, where
\[
\lb{v-}
 h={Q\/|k|_1}, \qq Q=\max\{\|q\|,\|q\|_1\},\qq
                                                |k|_1=\max\{1,|k|\};
\]
where $\|q\|=\int_0^\g |q(t)|dt$ and $\|q\|_1=\int_0^\g t|q(t)|dt$, and
\[
\begin{aligned}
\lb{y}
|y(x,k)|\le e^{2(\g-x)v_-+h},\\
|y(x,k)-1|\le h e^{2(\g-x)v_-+h},\\
|y(x,k)-1-y_1|\le {h^2\/2} e^{2(\g-x)v_-+h}.
\end{aligned}
\]
Moreover, $w\in \cE_\g$ and satisfies for any $k\in \C$:
\[
\lb{w1}
w(k)=i2k-q_0 +w_*(k), \qqq w_*(k)=-\int _0^\g q(t)\big(y(t,k)-1\big)dt,
\]
\[
\lb{w2}
|w_*(k)|\le \|q\| he^{h+2\g v_-},
\]
where $q_0=\int_0^\g q(t)dt$.

\end{lemma}
{\bf Proof.} The function $G(t,k)={\sin kt\/k}e^{ikt}$ satisfy
\[
\lb{y3}
|G(t,k)|\le {e^{2 v_-t}\/|k|_1}\ca 1 & if \ |k|\ge 1\\
                                   t & if \ |k|<1 \ac,
 \qqq t\ge 0, \ k\in \C,
\]
Consider the first case in \er{y3}: $|k|\ge 1$, the proof for the second case
$|k|<1$ is similar.
Substituting the estimate in \er{y3} for $|k|\ge 1$ into the identity
$$
y_n(x,k)=\int\limits_{x=t_0<t_1< t_2<...< t_n}
\lt(\prod\limits_{1\le j\le n}
G(t_{j}-t_{j-1},k)q(t_j)\rt)dt_1dt_2...dt_n,
$$
we obtain
\[
\begin{aligned}
|y_n(x,k)|\le {1\/|k|_1^n}\int\limits_{x=t_0<t_1< t_2<...< t_n}
\lt(\prod\limits_{1\le j\le n}e^{2v_-(t_{j}-t_{j-1})}
|q(t_j)|\rt)dt_1dt_2...dt_n\\
= {e^{-2v_- x}\/|k|_1^n}\int\limits_{x=t_0<t_1< t_2<...< t_n}
\lt(\prod\limits_{1\le j\le n}
|q(t_j)|\rt) e^{2v_-t_{n}}dt_1dt_2...dt_n\\
\le {e^{2(\g-x)v_-}\/|k|_1^n}\int\limits_{x=t_0<t_1< t_2<...< t_n}
|q(t_1)q(t_2)....q(t_n)| dt_1dt_2...dt_n\le e^{2(\g-x)v_-}{\|q\|^n\/n! |k|_1^n}.
\end{aligned}
\]

This shows that for each $x\ge 0$
the series \er{y2} converges uniformly on any bounded subset of $\C$.
Each term of this series is an entire function. Hence the sum is
an entire function. Summing the majorants we obtain estimates  \er{y}.
Thus the functions  $w, s$ and $\p _{+}(x,\cdot),\p _{+}'(x,\cdot),x\in\R$
are entire

The identity \er{3.1} gives
\[
\begin{aligned}
\lb{y5}
\p_+(0,k)'=ik-\int _0^\g\cos [kt] \ q(t)
\p_+(t,k)dt, \\
ik\p_+(0,k)=ik+i\int _0^\g\sin [kt] \ q(t)\p_+(t,k)dt,
\end{aligned}
\]
which yields
\[
\begin{aligned}
\lb{y6}
w(k)=\p_-(0,k )\p_+'(0,k )-\p_-'(0,k )\p_+(0,k )
=ik\p_+(0,k ) +\p_+'(0,k )\\
=2ik-q_0+w_*(k), \qqq w_*(k)=-\int _0^\g q(t)\big(y(t,k)-1\big)dt.
\end{aligned}
\]
Due to \er{y} we have the estimate
\[
\begin{aligned}
\lb{y7}
|w_*(k)|\le \int _0^\g|q(t)|h e^{2(\g-t)v_-+h}dt\le
h e^{2\g v_-+h}\int _0^\g|q(t)|dt=h \|q\|e^{2\g v_-+h},
\end{aligned}
\]
which is \er{w2}. Together with \er{ws} this implies that $w\in \cE_\g$.
\hfill \BBox

\

The identity \er{w1} and estimate \er{w2} give
 \[
\lb{w1x}
\begin{aligned}
w(k)=i2k-q_0 +w_*(k), \qq |w_*(k)|\le {Q^2\/|k|_1}e^{{Q\/|k|_1}+2\g v_-},\qq
 \forall \ k\in \C,\\
 where \ \ q_0=\int_0^\g q(t)dt,\qqq Q=\max\{\|q\|,\|q\|_1\}.
 \end{aligned}
\]

\

We show the relations between the Cases 1, 2 and 3. In order to show
this we use a standard trick. For $q\in L^1(\R_+)$ with $\supp q\ss
[0,\g]$ we define the Schr\"odinger operator $\wt
H$ acting in $L^2(\R )$ with an even compactly
supported potential $\wt q$ given by
\[
\lb{wtq} \wt H=-{d^2\/dx^2}+\wt q,\qqq \qqq \wt q(x)=q(|x|), x\in \R.
\]
 The Schr\"odinger equation
$
  -f''+\wt q(x)f=k^2f,\ \ \  \ \ \ k\in \C \sm \{0\},
$ has unique solutions $\wt\p_{\pm}(x,k)$  such that $\wt\p_+(x,k
)=e^{ixk },\ \ x\ge \g$ and $\wt\p_-(x,k)=e^{-ik x},\ \ x\leq -\g$.
Note that the symmetry of the potential $\wt q$ yields
$$
\wt\p_+(x,k)=\p_+(x,k)=\wt\p_-(-x,k) \qqq \forall x\in [0,\g].
$$
This implies that the Wronskian $\wt w(k)$ for the potential $\wt q$
satisfies
\[
\lb{tw} \wt w(k)=\{\wt\p_+(x,k),
\wt\p_-(x,k)\}|_{x=0}=2\p_+(0,k)\p_+'(0,k).
\]
where $\p_+(0,k)$ and $\p_+'(0,k)$ are the Jost functions for the
Cases 2 and 3 respectively with the potential $q$. This identity is
very useful. For example, if we have estimate \er{3} for the Case 1,
then  \er{tw} gives the estimate \er{3} for the Case 2 and 3.

We obtain the estimate of the Wronskian $\wt w(k)$.
We need to define the potential
$q_1=\wt q(x-\g)$, which has the support $\supp q_1\ss [0,2\g]$.
Using Lemma \ref{Tw} and \er{w1}, \er{w2} we obtain
\[
\lb{wtw}
\begin{aligned}
\wt w(k)=i2k-2q_0 +\wt w_*(k), \\
|\wt w_*(k)|\le \|q_1\| {\wt Q^2\/|k|_1}e^{4\g v_-+{\wt Q\/|k|_1}},
\end{aligned}
\]
where the corresponding constant $\wt Q$ is given by
\[
\begin{aligned}
\lb{wtQ}
\wt Q=\max \{\|q_1\|, \|q_1\|_1\}=2\|q\|\max \{1, \g\},
\end{aligned}
\]
since
$$
\|q_1\|=\int_0^{2\g} |q_1(x)|dx=2\|q\|,\qqq \|q_1\|_1=\int_0^{2\g} x|\wt q(x-\g)|dx=
\int_{-\g}^{\g} (t+\g)|\wt q(t)|dt=2\g \|q\|.
$$

We describe few facts about the resonances in the disk $\{|k|<r\}$.

\begin{proposition}
\lb{Tr}
 Let $q\in L^1(\R)$ and  $\supp q\ss [0,\g]$ for some $\g>0$.

 i)  Let $r>0$ and  let $q$ satisfy
\[
\lb{q0}
\|q\|(1+he^{2\g r+h})<2r, \qqq {\rm where} \ h={\max \{\|q\|,\|q\|_1\}\/\max\{1,r\}}.
\]
Then the function $w$ has only one simple zero in the disk $\{|k|<r\}$.

ii) Let $1\le r\le {1\/2\g}$ and let $2\|q\|\le r$.
Then the function $w$ has only one simple zero in the disk $\{|k|<r\}$.


\end{proposition}
\no {\bf Proof}. i) Let $w_0=i2k$ and $|k|=r$.  Then the estimates \er{w1},
\er{w2} imply
\[
|w(k)-w_0(k)|\le \|q\| (1+he^{2\g+h})=|w_0(k)|C_0,\qqq
C_0={\|q\|\/2r}(1+he^{2\g r+h}).
\]
Hence, if $C_0<1$, then by Rouche's theorem, $w(k)$ has only one simple zero,
as $w_0=2ik$ in the disk $\{|k|<r\}$.

ii) If we assume that $f(h)=h+h^2e^{1+h}<2$, then \er{q0} holds true.

We have that the function $f$ is increasing. We get
the estimates $f({1\/2})<2$, since  $h={\|q\|\/r}\le {1\/2}$, which yield that $w$
has only one simple zero in the disk $\{|k|<r\}$.
\hfill\BBox

{\bf Remark.} 1) This proposition shows that for  any $\g, r>0$ and sufficiently
small $\|q\|$  the function $w$ has only one simple zero in the disk $\{|k|<r\}$.

2)  The number of zeros of $w$ in the disk with large $r$  depends on the diameter
$\g$. But it is very interesting to determine the asymptotics of $\cN(r,w)$ as
$r\to \iy$. This problem remains open for  long period.

3) We need to say that \er{q2} holds also true for the Case 2 and 3.
 Then \er{3} for Case 1 yields \er{3} for Cases 2 and 3.

\

 {\bf Proof of Theorem \ref{TrJS}}.
Lemma \ref{Tw} gives that the function $w\in \cE_\g$.

We apply the estimate \er{L1} to the function $w$. Thus  due to the estimate \er{w1x} and the simple fact $|q_0|\le Q$,
we obtain \er{q2}.

We show \er{q3}.
The identity \er{tw} gives $\cN (r,\p_+(0,\cdot))+\cN
(r,\p_+'(0,\cdot))=\cN (r,\wt w)$, where $\wt w$ is the Wronskian
for the potential $\wt q$ given by \er{wtq} and $\wt w$ satisfies \er{wtw}, \er{wtQ}. Then these arguments and  the estimate \er{q2} implies
 \er{q3}. \hfill $\BBox$

 \medskip

{\bf Proof of Theorem \ref{T1}.} We consider Case 1 and let $\supp
q\ss [0,\g]$. By Lemma \ref{Tw}, the function $w$ belongs to
$\cE_\g$ with $Q=\max\{\|q\|,\|q\|_1\}$ and satisfies \er{w1x}.
 Thus  the estimate \er{L3}, \er{L2}  from  Theorem
\ref{TEf} give \er{3} for  Case 1.

Consider Cases 2 and 3.  Then we have the identity $\wt w(k)=2\p_+(0,k)\p_+'(0,k)$,
where $\wt w$ is the Wronskian for the potential $\wt q$ given by \er{wtq}
(see  \er{tw}). Then  \er{3} for $\wt q$ (Case 1) together with  \er{wtQ} implies
\er{3} for Cases 2 and 3.
\hfill $\BBox$

\begin{lemma}
\lb{Tp} The function $Y_p=\int_\R {(1+x^2)^{-{p\/2}}}dx, p>1$ is
strongly monotonic and convex on $(1,\iy)$, since $Y_p'<0, Y_p''>0$,
and  satisfies \er{Yp}  and the identity \er{p1}, i.e., $Y_p=\sqrt \pi{\G({p-1\/2})/\G({p\/2})}$ for all $p>1$.
\end{lemma}
\no {\bf Proof.} We rewrite $Y_p$ in the form $Y_p=\int_\R e^{-p
f(x)}dx$, where $f(x)={1\/2}\log (1+x^2)$. This yields
$Y_p'=-\int_\R f(x)e^{-p f(x)}dx<0$ and $Y_p''=\int_\R f^2(x)e^{-p
f(x)}dx>0$. We rewrite $Y_p$ in another form
$$
Y_p=2\int_0^\iy {(1+x^2)^{-{p\/2}}}dx=\int_0^\iy
{t^{-{1\/2}}(1+t)^{-{p\/2}}}dt=\sqrt \pi{ \G(z)\/\G(z+{1\/2})}, \qq
z={p-1\/2},
$$
where $\G$ is the Gamma function, see identities (2), (5) in Section
1.5 in \cite{EMOT53}.

Asymptotics  ${ \G(z)\/\G(z+{1\/2})}=z^{-{1\/2}}(1+O({1\/z}))$ as
$z\to \iy$, see  (5) in Sect. 1.18 in \cite{EMOT53} give \er{Yp} as
$p\to \iy$. Recall that $\G$ is meromorphic in $\C$ with simple
poles $0,-1,-2,..$. Then the identities $\G(1+z)=z\G(z)$ and
$\G(1)=1, \G(1/2)=\sqrt \pi$ yield  \er{Yp} as $p\to 1$:
$$
Y_p=\sqrt \pi{ \G(z)\/\G(z+{1\/2})}={\sqrt\pi \G(1+z)\/
z\G(z+{1\/2})}={1+O(z)\/z}\qq \as \  \ z={p-1\/2}\to 0,
$$
 \hfill $\BBox$

 \medskip

\no {\bf Acknowledgments.} \footnotesize I am grateful to  Michael
Loss, Atlanta, Grigori Rozenblum, G\"oteborgs, and Sergey Morozov,
Munich, for useful comments about Lieb-Thirring inequalities and
Elliot Lieb, Princeton, for the reference \cite{DFLP06}. Moreover, I
am also grateful to  Alexei Alexandrov  (St. Petersburg) for
stimulating discussions about the Carleson Theorem and about the
absolute constant $C$ in \er{Ce}.

This work was supported by the Ministry of education and science of
the Russian Federation, state contract 14.740.11.0581 and the RFFI
grant "Spectral and asymptotic methods for studying of the
differential operators" No 11-01-00458.


\begin{thebibliography}{9}
\footnotesize


\bibitem   {C58} Carleson, L. An interpolation problem for bounded
 analytic functions.  Amer. J. Math. 80 (1958), 921--930.

\bibitem  {C62} Carleson, L.
  Interpolations by bounded analytic functions and the corona problem.
   Ann. of Math.  76 (1962), 547--559.


\bibitem{BKW03}  Brown, B.;  Knowles, I.; Weikard, R.
On the inverse resonance problem, J. London Math. Soc.  68 (2003),
no. 2, 383--401.


\bibitem  {G65} S. Golden, Lower bounds for the Helmholtz function,
Phys. Rev. (2) 137 (1965), B1127--B1128.

\bibitem   {DT79} Deift, P.; Trubowitz, E. Inverse scattering on
the line, Commun. Pure and Applied Math., 32(1979), 121--251.



\bibitem  {DFLP06} Dolbeault, J.; Felmer, P.; Loss, M.; Paturel,
E. Lieb-Thirring type inequalities and Gagliardo-Nirenberg inequalities
 for systems. J. Funct. Anal. 238 (2006), no. 1, 193--220.

\bibitem  {EMOT53} Erdelyi, A.; Magnus, W.; Oberhettinger, F.;
Tricomi, F. Higher transcendental functions. Vol. I, Based, in part,
on notes left by Harry Bateman.
 McGraw-Hill Book Company, Inc., New York-Toronto-London, 1953.


\bibitem   {H99} Hitrik, M. Bounds on scattering poles in one dimension,
Commun. Math. Phys. 208(1999), 381--411.

\bibitem   {F97} Froese, R. Asymptotic distribution of resonances in one
dimension. J. Diff. Eq. 137 (1997), no. 2, 251--272.

\bibitem     {G81}  Garnett, J. Bounded analytic functions, Academic Press,
New York, London, 1981.

\bibitem   {J64}   Javrjan, V. A. Some perturbations of self-adjoint operators.
  (Russian) Akad. Nauk Armjan. SSR Dokl. 38 (1964), 3--7.


\bibitem   {JP51} Jost, R.; Pais, A. On the scattering of a particle
 by a static potential.
 Phys. Rev., 82 (1951). 840--851.

\bibitem   {Koo88} Koosis, P. The logarithmic integral I,
Cambridge Univ. Press, Cambridge, London, New York 1988.

\bibitem {K04} Korotyaev, E. Inverse resonance scattering  on the half
line, Asymptotic Anal.  37(2004), No 3/4, 215--226.


\bibitem   {K04s}  Korotyaev, E. Stability for inverse resonance problem.
 Int. Math. Res. Not. 2004, no. 73, 3927--3936.


\bibitem   {K05}  Korotyaev, E. Inverse resonance scattering on the
real line. Inverse Problems 21 (2005), no. 1, 325--341.

\bibitem  {K11} Korotyaev, E. Resonance theory for perturbed Hill
operator, Asymp. Anal. 74(2011), No 3-4, 199--227.

\bibitem{L93} Levin, B. Ya. Lectures on entire functions.
  Translations of Mathematical Monographs, 150.
  American Mathematical Society, Providence, RI, 1996.




\bibitem {LT76}  Lieb E., Thirring, W. Inequalities for the moments of
the eigenvalues of the Schr\"odinger
Hamiltonian and their relation to Sobolev inequalities. Studies in
Math. Phys., Essays in Honor of Valentine Bargmann, Princeton 1976.

\bibitem   {L13} Loss, M. private communication.


\bibitem   {Ma86} Marchenko, V. Sturm-Liouville operator and applications.
Basel, Birkh\"auser 1986.

\bibitem   {MSW10}
 Marletta, M.;  Shterenberg, R.; Weikard, R.,
On the Inverse Resonance Problem for Schr\"odinger Operators,
Commun. Math. Phys., 295(2010), 465--484.


\bibitem   {Me85} Melin, A. Operator methods for inverse scattering on the
 real line. Comm. P.D.E. 10(1985), 677--786.

\bibitem  {M83}  Melrose, R. Polynomial bound on the number of scattering
 poles, J. Funct. Anal. 53(1983), 287--303.

\bibitem   {S00} Simon, B. Resonances in one dimension and Fredholm
determinants,  J. Funct. Anal. 178 (2000), no. 2, 396--420.



 \bibitem {SZ91} Sj\"oostrand, J.; Zworski, M. Complex scaling and the
 distribution of scattering poles. J. Amer. Math. Soc. 4 (1991), no. 4, 729-769.

\bibitem  {T65}  Thompson, C. Inequality with applications in statistical
 mechanics, J. Math. Phys. 6(1965), 1812--1813.




\bibitem  {Z87} Zworski, M. Distribution of poles for scattering
on the real line, J. Funct. Anal. 73(1987), 277--296.

\bibitem {Z89}  Zworski, M. Sharp polynomial bounds on the number
of scattering poles of radial potentials. J. Funct. Anal.  82
(1989), no. 2, 370--403.


\bibitem {Z02} Zworski, M. SIAM, J.  Math. Analysis, "A remark on
isopolar potentials" 82(2002), No 6,  1823--1826.



\end{thebibliography}
\end{document}